\documentclass[a4paper, 12pt]{article}
\usepackage{enumerate, theorem}
\usepackage{amsmath, amsfonts, amssymb}
\usepackage[height=22.5cm, width=15cm]{geometry}
\usepackage[all]{xy}
\usepackage{mathrsfs, yfonts}

\DeclareMathOperator{\Sym}{Sym}

\DeclareMathOperator{\C}{\mathbb{C}}

\DeclareMathOperator{\f}{\mathcal{C}_n^f(M)}

\DeclareMathOperator{\M'}{\mathcal{C}_n(M,M')}

\newcommand{\parag}[1]{\paragraph{\sc{#1.}}}

\newtheorem{thm}{Theorem}[subsection]
\newtheorem{defn}[thm]{Definition}
\newtheorem{cor}[thm]{Corollary}
\newtheorem{prop}[thm]{Proposition}
\newtheorem{lemma}[thm]{Lemma}

\begin{document}

\title{Quasi-proper meromorphic equivalence relations.}

\date{31/05/10.}

\author{Daniel Barlet\footnote{Barlet Daniel, Institut Elie Cartan UMR 7502  \newline
Nancy-Universit\'e, CNRS, INRIA  et  Institut Universitaire de France, \newline
BP 239 - F - 54506 Vandoeuvre-l\`es-Nancy Cedex.France. r\newline
e-mail : barlet@iecn.u-nancy.fr}.}

\maketitle

 \hfill {\it En hommage \`a H. Grauert et R. Remmert.}
 
 \parag{Abstract}
 The aim of this article is to complete results of [M.00] and  [B.08] and to show that they  imply a rather general existence theorem for meromorphic quotient of strongly quasi-proper meromorphic equivalence relations. In this context, generic equivalence classes are asked to be pure dimensionnal closed analytic subset with finitely many irreducible components. As an application of these methods we prove a Stein factorization theorem for a strongly quasi-proper map.
 
 \bigskip
 
 \parag{AMS Classification} 32 C 25, 32 H 35, 32 H 04 .
 \parag{Key words} Meromorphic equivalence relations, meromorphic quotients, geometric f-flattening, (strongly) quasi-proper.
 
 \newpage
 
 \tableofcontents
 
 \section{Introduction.}
 
 The aim of this article is to complete results of [M.00] and  [B.08] and to show that they  imply a rather general existence theorem for meromorphic quotient of quasi-proper meromorphic equivalence relations. We try also to put some more light on the "topological" condition which is needed to have such a meromorphic quotient. I hope that these stronger results and this new presentation of this topological condition will help potential users.\\
 Of course the results here are in some sense a wide generalization of the classical Henri Cartan's quotient theorem [C.60]. But without requiring compactness of equivalence classes and in a more "geometrical" spirit which allows for instance, to always have a meromorphic quotient in the proper case. Remember that to have a functorial (holomorphic) quotient as complex space, H. Cartan gaves a necessary and sufficient condition which is not automatic.\\
 The method we use is deeply related to the study of geometric f-flattening of a quasi-proper map, which is a generalization in the case of non compact fibers of the geometric flattening theorem of [B.78]  for proper holomorphic maps.  We can only prove the existence of a meromorphic quotient with a strong quasi-properness assumption, and in this case  the corresponding meromorphic quotient map admits a geometric f-flattening. It is easy to see that there exists examples where a quasi-proper meromorphic quotient exists such that  the quotient map does not admit geometric f-flattening. But this phenomenon is related to the fact that quasi-properness is a notion which is too weak for non equidimensional map : \\
 the fact that all irreducible components of a big fiber meet a compact set does not imply that the generic nearby fibers have not irreducible components which  escape to infinity (see [M.00] or [B.08] for the phenomenon of "escape at infinity"). \\
 The topological condition we add garanties that this pathology does not happen and gives in fact a notion of "strong quasi-properness" which is equivalent to the (local)  existence of geometric f-flattening. This corresponds also to the fact that the generic fibers of the map may be completed in  a f-meromorphic family of cycles. As quasi-properness this notion is local on the target space, but, in contrast with the  quasi-proper case, this condition is stable by proper modification of the image.\\
 Another key point for proving this rather general existence theorem is the reparametrisation theorem of D. Mathieu [M.00] which is a  consequence of his generalization of N. Kuhlmann's semi-proper direct image theorem [K.64] and [K.66]  in the case where the target space is an open set of a Banach space. We show in the appendix how it implies in fact a semi-proper direct image theorem\footnote{but of course the functor "f-analytic families of \ $n-$cycles in \ $M$" is not representable in the category of reduced complex spaces.} with values in \ $\mathcal{C}_n^f(M)$ \ the space of finite type closed \ $n-$cycles in a complex space \ $M$.\\
 
 As an application of these ideas, we give a Stein's factorization theorem for a strongly quasi-proper holomorphic map.\\
 
 To conclude this introduction let me say that I consider  this work as a tribute to professors R. Remmert and H. Grauert :   the direct image theorem of R. Remmert [R.57]  is the basic idea used here to produce meromorphic quotients and  I think that the present work is also a far reaching conclusion on this problem which was initiated by H. Grauert [G.83] and  [G.86] and his student B. Siebert [S.93] and [S.94].\\
 
 During the final draft of this article the first  discussions about the article [B.Mg.10] were going on and several interesting points were clarified. So I want to thanks J\'on Magn\'usson for his help.

\section{Quasi-proper geometric flattening.}

\subsection{Geometric flattening.}

Let \ $f : M \to S$ \ be a holomorphic map between two reduced and irreducible complex spaces and put \ $n : = \dim M - \dim S$.\\
 Recall that such a map is  called {\bf geometrically flat} when there exists an analytic family \ $(X_s)_{s \in S}$ \ of \ $n-$cycles in \ $M$ \ parametrized by \ $S$ \ such that for any \ $s \in S$ \ the support \ $\vert X_s\vert$ \ of the cycle \ $X_s$ \ is the set theoretic fiber \ $f^{-1}(s)$ \ and such that, for general \ $s$,  the cycle \ $X_s$ \ is reduced (i.e. equal to its support, so each irreducible component is of multiplicity 1 in this cycle). Notice that such a map \ $f$ \  is open and that, in the case where \ $S$ \ is normal,  the map \ $f$ \ is geometrically flat if and only if it is \ $n-$equidimensional (see [B.75] or [BOOK]).\\
For a geometrically flat \ $f$ \ we have a classifying map for the analytic family \ $(X_s)_{s \in S}$ \ which is a map \ $ \varphi : S \to \mathcal{C}_n^{loc}(M)$ \ given by \ $\varphi(s) : = X_s$.

\begin{defn}\label{geom.flattning}
Given an arbitrary holomorphic map \ $f : M \to S$ \ between two reduced and irreducible complex spaces, a {\bf geometric flattening} for \ $f$ \ is a  proper holomorphic modification \ $\tau : \tilde{S} \to S$ \ of \ $S$ \ such that the strict transform\footnote{In general \ $\tilde{M}$ \ is, by definition, the union of irreducible components of \ $\tilde{S}\times_S M$ \ which dominate an irreducible component of \ $S$. When  \ $M$ \ is irreducible and \ $f$ \ surjective \ $\tilde{M}$ \ is the irreducible component of \ $\tilde{S}\times_S M$ \ which  surjects on \ $S$.}  \\ $\tilde{f} : \tilde{M} \to \tilde{S} $ \ is a geometrically flat holomorphic map.
\end{defn}

As the holomorphic map \ $\pi : \tilde{M} \to M$ \ is proper we have (see [B.75] or [BOOK]) a direct image map \ $\pi_* : \mathcal{C}_n^{loc}(\tilde{M}) \to \mathcal{C}_n^{loc}(M)$ \ which is holomorphic in the sense that it preserves the analyticity of families of \ $n-$cycles. So when we have a geometric flattening for \ $f$ \ the classifying map \ $\tilde{\varphi} : \tilde{S} \to \mathcal{C}_n^{loc}(\tilde{M})$ \ composed with the direct image \ $\pi_*$ \ gives a map \ $\pi_*\circ \tilde{\varphi} : \tilde{S} \to \mathcal{C}_n^{loc}(M)$ \ which will classify the "fibers" of \ $f$. The existence of a geometric flattening for the map \ $f$ \ may be considered as the meromorphy of the classifying map of the analytic family of generic fibers of \ $f$ \ along the center \ $\Sigma$ \ of the proper holomorphic  modification \ $\tau$ \ (see the definition \ref{f-mero.} below).\\

\subsection{Quasi-proper geometric flattening.}

In what follows we shall consider mainly finite type \ $n-$cycles, that is to say \ $n-$cycles having finitely many irreducible components. Recall that the classical corresponding  relative notion is given by the following definition

\begin{defn}\label{Quasi-proper}
Let \ $f : M \to S$ \ be a holomorphic map between two reduced complex spaces. We say that \ $f$ \ is {\bf quasi-proper} if for any point \ $s_0 \in S$ \ there exists an open neighbourhood \ $S'$ \ of \ $s_0$ \ in \ $S$ \ and a compact set \ $K$ \ in \ $M$ \ such that for any \ $s \in S'$ \ any irreducible component \ $\Gamma$ \ of \ $f^{-1}(s)$ \ meets \ $K$.
\end{defn}

Of course the fibers of a quasi-proper map are finite type cycles and with some "uniform" local condition for the finiteness of the irreducible components of the fibers. The notion of quasi-proper map is not topological. Nevertheless, it may be defined for a continuous map \ $f : M \to S$ \ where \ $M$ is a complex space and \ $S$ \ a topological space, if we know that any fiber of \ $f$ \ is an analytic subset of \ $M$. This is, for instance, the case for the projection on \ $S$ \ of  the (set theoretic) graph \ $\vert G \vert \subset S \times M$ \ of a continuous family of \ $n-$cycles \ $(X_s)_{s \in S}$ \ of \ $M$ \ parametrized by a Hausdorff topological space \ $S$.\\

The following notion is purely topological.

\begin{defn}\label{semi-proper}
Let \ $f : M \to S$ \ be a continuous map between two Hausdorff topological spaces. We say that \ $f$ \ is  {\bf semi-proper} if for any point \ $s_0 \in S$ \ there exists an open neighbourhood \ $S'$ \ of \ $s_0$ \ in \ $S$ \ and a compact set \ $K$ \ in \ $M$ \ such that \ $f(M)\cap S' = f(K) \cap S'$.
\end{defn}

Of course a quasi-proper map is always semi-proper. Kuhlmann's theorem (see [K.64] and [K.66] ) generalizes Remmert 's direct image theorem [R. 57] to the semi-proper case :

\begin{thm}[Kuhlmann]
Let \ $f : M \to S$ \ be an holomorphic map between two reduced complex spaces. Assume that \ $f$ \ is semi-proper. Then \ $f(M)$ \ is a (closed)  analytic subset of \ $S$.
\end{thm}

Recall now that we introduced  in [B.08], for any given complex space \ $M$, the topological space \ $\mathcal{C}_n^f(M)$ \ of finite type \ $n-$cycles with a topology finer that the topology induced by the obvious inclusion \ $\mathcal{C}_n^f(M) \subset \mathcal{C}_n^{loc}(M)$, where \ $\mathcal{C}_n^{loc}(M)$ \ is the (topological) space of closed $n-$cycles in \ $M$. We defined also
 the notion of an f-analytic family of (finite type) \ $n-$cycles in \ $M$. Let \ $S$ \ be a reduced complex space. The definition of the topology on \ $\mathcal{C}_n^f(M)$ \ is given in order to have  that  f-analytic families of finite type \ $n-$cycles in \ $M$ \ are exactely the analytic families \ $(X_s)_{s \in S}$ \  of \ $n-$cycles in \ $M$ \ that satisfy the following condition :
\begin{itemize}
\item The projection of the support \ $\vert G\vert \subset S \times M$ \ of the graph \ $G$ \ of the family is quasi-proper on \ $S$.
\end{itemize}
It of course implies that for each \ $s \in S$ \ the \ $n-$cycle \ $X_s$ \ is of finite type (i.e. has finitely many irreducible components). \\
This is a purely topological requirement on the family, corresponding to the continuity of the family for the finer topology defined on \ $\mathcal{C}_n^{f}(M)$.\\
For a f-analytic family the classifying map \ $\varphi : S \to \mathcal{C}_n^{loc}(M)$ \ factors through a continuous classifying map \ $\varphi^f : S \to \mathcal{C}_n^f(M)$ \ and the (continuous)  inclusion 
 $$ i : \mathcal{C}_n^f(M) \to \mathcal{C}_n^{loc}(M) .$$
 
 \begin{defn}\label{geometrically f-flat}
 Let \ $M$ \ and \ $S$ \ be two reduced complex spaces and let \ $f : M \to S$ \ be a quasi-proper holomorphic map which is \ $n-$equidimensional. We shall say that \ $f$ \ is { \bf geometrically f-flat} when there exists an f-analytic family \ $(X_s)_{s \in S}$ \ of \ $n-$cycles in \ $M$ parametrised by \ $S$ \ such that we have for each \ $s \in S$ \ the equality \ $\vert X_s\vert = f^{-1}(s)$ \ and such that for generic \ $s \in S$ \ the cycle \ $X_s$ \ is reduced (so each irreducible component is of multiplicity 1 in \ $X_s$).
 \end{defn}
 
 It is easy to see, using local compactness of \ $S$ \  that a geometrically flat map \ $f : M \to S$ \ is geometrically f-flat if and only if \ $f$ \ is quasi-proper.\\
 
 For a geometrically f-flat map \ $f$ \ we have an "holomorphic"  classifying map
 $$ \varphi^f : S \to \mathcal{C}_n^f(M) $$
 associated to the f-analytic family of "fibers" of \ $f$.
 
 \bigskip

Let us now consider a quasi-proper surjective holomorphic map \ $f : M \to S$ \ between two reduced and irreducible complex spaces, and let \ $n : = \dim M - \dim S$. As the set of points \ $x$ \  in \ $M$ \ where the fiber at \ $x$ \  has dimension \ $> n$ \ is a closed analytic subset of \ $M$ \ which is \ $f-$saturated in the sense that, in each fiber of \ $f$, this subset  is an union of irreducible components of the fiber, its image is an analytic subset \ $\Sigma$ \ in \ $S$ \ which has no interior point in \ $S$. This is a consequence of Kuhlmann's theorem, because the restriction of  a  quasi-proper map \ $f$ \ to a \ $f-$saturated analytic subset is again quasi-proper.\\
Without loss of generality, we may assume that the non normal points in \ $S$ \ are in \ $\Sigma$. And now the restriction of \ $f$ \ to \ $M \setminus f^{-1}(\Sigma)$ \ is quasi-proper and equidimensional on the normal complex space \ $S \setminus \Sigma$. So we have a f-analytic family of \ $n-$cycles \ $(X_s)_{s \in S \setminus \Sigma}$ \ associated to the fibers of this map and a corresponding classifying map
$$ \varphi^f : S\setminus \Sigma \to \mathcal{C}_n^f(M). $$

Recall that a subset \ $Q $ \ in \ $\mathcal{C}_n^f(M)$ \ is compact if and only if it is compact  in \ $\mathcal{C}_n^{loc}(M)$ \ and such that the topologies induced by \ $\mathcal{C}_n^{loc}(M)$ \ and\ $\mathcal{C}_n^f(M)$ \ co{\"i}ncide. In general the uniform  boundness of the volume in each compact set of \ $M$, which is equivalent, thanks to Bishop's theorem [Bi.64], to the relative compactness in \ $\mathcal{C}_n^{loc}(M)$ \ (see [BOOK] ch.IV),  will be easy to check. But to check whether the topologies induced by \ $\mathcal{C}_n^f(M)$ \  and \ $\mathcal{C}_n^{loc}(M)$ \ co{\"i}ncide (which is equivalent to the "non escape at infinity") is  a sticky point in the sequel. For a precise description of the topologies on \ $\mathcal{C}_n^{loc}(M)$ \ and \ $\mathcal{C}_n^f(M)$ \ and their comparison see [B.08].\\

The following lemma gives a precise characterization of a compact subset in \ $\mathcal{C}_n^f(M)$ \ without any reference to the "escape at infinity". 

\begin{lemma}\label{compact in f-}
A closed subset \ $\mathcal{B}$ \ in \ $\mathcal{C}_n^f(M)$ \ is  compact in \ $\mathcal{C}_n^f(M)$ \ if and only if it is a compact subset of \ $\mathcal{C}_n^{loc}(M)$ \ and there exists a compact set \ $K \subset M$ \ such that any irreducible component of any cycle in \ $\mathcal{B}$ \ meets \ $K$.
\end{lemma}

\parag{Proof}
A compact subset in \ $\mathcal{C}_n^f(M)$ \ is clearly  a compact in \ $\mathcal{C}_n^{loc}(M)$ \ because the inclusion map is continuous. For any point \ $X \in \mathcal{B}$ \ let \ $W_X$ \ be a relatively compact open set  in \ $M$ \ which meets all irreducible components of \ $X$. Then \ $\Omega(W_X)$, the set of cycles in \ $\mathcal{C}_n^f(M)$ \ such that each  irreducible component meets \ $W$,  is an open neighbourhood of \ $X$ \ in \ $\mathcal{C}_n^f(M)$. Choosing a finite sub-covering \ $(\Omega(W_{X_i}))_{i \in I}$ \ of the covering of \ $\mathcal{B}$ \ by the open sets \ $\Omega(W_X)_{X \in \mathcal{B}}$, gives a relatively compact open set \ $W : = \cup_{i \in I} \ W_{X_i}$ \ such that any irreducible component of any \ $Y \in \mathcal{B}$ \ meets \ $W$. \\
Conversely, if a closed subset \ $\mathcal{B} $ \   in \ $\mathcal{C}_n^f(M)$ \ is compact in \ $\mathcal{C}_n^{loc}(M)$ \   and if any irreducible component of any cycle in \ $\mathcal{B}$ \ meets a compact \ $K$ \ in \ $M$, to show   compactness in \ $\mathcal{C}_n^f(M)$ \ take any sequence \ $(X_{\nu})_{\nu \in \mathbb{N}}$ \ in \ $\mathcal{B}$. Up to pass to a subsequence, we may assume that \ $(X_{\nu})_{\nu \in \mathbb{N}}$ \ converges to a cycle \ $X \in \mathcal{B}$ \ in the topology of \ $\mathcal{C}_n^{loc}(M)$. We want to show that the convergence takes place in the sense of the topology of \ $\mathcal{C}_n^f(M)$.\\
If \ $X$ \ is the empty \ $n-$cycle, cover \ $K$ \ with finitely many \ $n-$scales (always adapted to \ $\emptyset$). Then for \ $\nu \gg 1$ \ the degree of \ $X_{\nu}$ \ in each of these scales has to be \ $0$ \ and so \ $\vert X_{\nu}\vert \cap K = \emptyset$. So the only possibility is that \ $X_{\nu} = X$ \ for \ $\nu \gg 1$.\\
 When \ $X$ \ is not the empty \ $n-$cycle, we have to prove that if any irreducible component of \ $X$ \ meets an open set \ $W$, then for \ $\nu \gg 1$ \ any irreducible component of \ $X_{\nu}$ \ also meets \ $W$. But if, for an infinite sequence of \ $\nu \geq 0$, there exists an irreducible component \ $C_{\nu}$ \ of \  $X_{\nu}$ \ which does not meet \ $W$, choose a point \ $x_{\nu} \in K \cap C_{\nu}$. The points \ $x_{\nu}$ \ are then in \ $K \setminus W$ \ which is compact. So, up to pass to a subsequence we may assume, first that the sequence \ $(C_{\nu})_{\nu \in \mathbb{N}}$ \ converges to a cycle \ $Y$ \ in the sense of the topology of \ $\mathcal{C}_n^{loc}(M)$ \ and that  the sequence\ $(x_{\nu})_{\nu \in \mathbb{N}}$ \ converges to \ $x \in K \setminus W$. But then we have \ $\vert Y \vert \subset \vert X\vert $ \ and \ $x \in \vert Y\vert$. So \ $Y$ \ is not the empty cycle and any of its irreducible components meet \ $W$. So for \ $\nu$ \ large enough \ $C_{\nu}$ \ meets \ $W$. Contradiction. As \ $\mathcal{B}$ \ is closed in \ $\mathcal{C}_n^f(M)$ \ we conclude that the sequence \ $(X_{\nu})_{\nu \in \mathbb{N}}$ \ converges to \ $X \in \mathcal{B}$ \ in the topology of \ $\mathcal{C}_n^f(M)$.$\hfill \blacksquare$\\

\parag{Remark} If \ $\mathcal{B}$ \ is a compact set in \ $\mathcal{C}_n^f(M)$ \ then the subset $$\widehat{\mathcal{B}} : = \{ X \in \mathcal{C}_n^f(M) \ / \  \exists Y \in \mathcal{B} \quad  X \leq Y \}$$
is also compact. The compactness of \ $\widehat{\mathcal{B}}$ \  in \ $\mathcal{C}_n^{loc}(M)$ \ is obvious, the same compact \ $K \subset M$ \ which meets all irreducible components of a cycle in  \ $\mathcal{B}$ \ meets also each irreducible component of a cycle in \ $\widehat{\mathcal{B}}$; finally this subset is closed in \ $\mathcal{C}_n^f(M)$ \ thanks to the lemma \ref{limite inegalite}.\\
Note that \ $\emptyset$ \ is always in \ $\widehat{\mathcal{B}}$ \ but that \ $\widehat{\mathcal{B}} \setminus \{\emptyset\}$ \ is closed in \ $\widehat{\mathcal{B}}$ \ as \ $\{\emptyset\}$ \ is open in \ $\mathcal{C}_n^f(M)$ ;  so \ $\widehat{\mathcal{B}} \setminus \{\emptyset\}$ \ is also a compact set.\\

Now comes the main difference between the fact that we consider arbitrary closed \ $n-$cycles or finite type closed \ $n-$cycles. To make this difference transparent, let me use the following definition :

\begin{defn}\label{properly extendable map}
Let \ $S$ \ be an irreducible complex space and \ $\Sigma \subset S$ \ be a closed analytic subset with no interior point in \ $S$. Let \ $\varphi : S \setminus \Sigma \to Z$ \ be a continuous map to a Hausdorff topological space \ $Z$. We say that the map \ $\varphi$ \ is {\bf properly extendable along} \ $\Sigma$ \ if there exists a Hausdorff topological space \ $Y$, a continuous  map \ $\sigma : Y \to S$ \ and a continuous map
$$ \psi : Y \to Z$$
such that the map \ $\sigma$ \ is a { \bf proper topological modification along \ $\Sigma$} \ and such that  \ $\psi$ \ extends \ $\varphi$. Topological modification means  that \ $\sigma$ \  is continuous and  proper, that the set \ $\sigma^{-1}(\Sigma)$ \ has no interior point in \ $Y$ \ and that the restriction
$$ \sigma : Y \setminus \sigma^{-1}(\Sigma) \to S \setminus \Sigma$$
is a homeomorphism. The fact that \ $\psi$ \ extends \ $\varphi$ \ means that on \ $S \setminus \Sigma$ \ we have \ $\varphi = \psi\circ \sigma^{-1} $.
\end{defn}

Now we have the following key theorem :

\begin{thm}\label{flattening 1} {\bf [Geometric f-flattening; first version.]}
Let \ $M$ \ and \ $S$ \ be two reduced and irreducible complex spaces. Let \ $f : M \to S$ \  be a holomorphic map and assume that \ $\Sigma \subset S$ \ is a (closed) analytic subset in \ $S$ \ with no interior point, containing the non normal points of \ $S$ \ and such that the restriction of \ $f$ \ to \ $M \setminus f^{-1}(\Sigma)$ \ is \ $n-$equidimensional, where \ $n : = \dim M - \dim S$.
\begin{enumerate}[i)]
\item The map \ $\varphi : S \setminus \Sigma \to \mathcal{C}_n^{loc}(M) $ \ classifying the fibers of \ $f$ \ over \ $S \setminus \Sigma$ \  is always  properly extendable along \ $\Sigma$.\\
\item Assume that \ $f : M \to S$ \ is quasi-proper, and let \ $\varphi^f : S \setminus \Sigma \to \mathcal{C}_n^f(M)$ \ be the classifying map of the fibers of \ $f$ \ over \ $S \setminus \Sigma$. \\
If \ $\varphi^f$ \ is properly extendable along \ $\Sigma$ \ then there exists a proper holomorphic modification \ $\tau :\tilde{S} \to S$ \ and an f-analytic family of finite type \ $n-$cycles in \ $M$ \ which extends \ $\varphi^f$ \ to \ $\tilde{S}$.\\
\end{enumerate}
\end{thm}

 So in the case of arbitrary closed \ $n-$cycles, the continuous extension of the classifying map \ $\varphi$ \  to a proper topological modification of \ $S$ \ along \ $\Sigma$ \ is always possible, but this does not allow us to obtain a holomorphic extension on a proper holomorphic modification of \ $S$ along \ $\Sigma$. \\
In the case of finite type \ $n-$cycles the quasi-properness asumption on \ $f$ \ does not imply automatic topological extension for the f-classifying map \ $\varphi^f$. But when this topological extension exists, the f-classifying map  can be holomorphically extended on a suitable proper holomorphic modification of \ $S$ \ along \ $\Sigma$.

\parag{Proof of the theorem} Let us begin by the case i): define \ $\Gamma \subset S \times \mathcal{C}_n^{loc}(M)$ \ as the closure of the graph of \ $\varphi$. Then \ $\Gamma$ \ is \ $S-$proper: this is an easy consequence of the characterisation of compact sets in \ $\mathcal{C}_n^{loc}(M)$ \ via Bishop's theorem (see [Bi.64] and  [BOOK] ch.IV), and the result of [B.78] on the local boundness of the volume of generic fibers for an holomorphic map ;  so we want to prove two facts\footnote{Recall that we don't know that  \ $\Gamma$ \ is locally compact ; so \ $p$ \ proper means that \ $p$ \ is  closed with compact fibers.} 
\begin{enumerate}[1)]
\item The projection \ $p : \Gamma \to S$ \ is a closed map ;
\item its fibers are compact subset in \ $\mathcal{C}_n^{loc}(M)$.
\end{enumerate}
Let \ $F$ \ be a closed set in \ $\Gamma$, and assume that the sequence \ $(s_{\nu})_{\nu \in \mathbb{N}}$ \ in \ $p(F)$ \ converges to \ $\sigma \in S$. Let \ $(s_{\nu}, C_{\nu}) \in F$ \ for \ $\nu \in \mathbb{N}$. Now fix a compact neighbourhood \ $L$ \ of \ $\sigma$ \ in \ $S$ ;  then using [B.78], for each compact set \ $K \subset M$ \ and any fixed continuous hermitian metric \ $h$ \ in \ $M$ \ we may find  a constant \ $\gamma(K,h)$ \   such for any \ $s \in L \setminus \Sigma$ \  we have \ $vol_h(K \cap \varphi(s)) \leq \gamma(K,h)$. This inequality extends by continuity to the closure of \ $\varphi(L \setminus \Sigma)$ \ in \ $\mathcal{C}_n^{loc}(M)$, so to  \ $p^{-1}(L)$ \ in \ $F$. This implies that \ $p^{-1}(L)$ \ is a compact set in \ $F$. This allows us, up to pass to a subsequence, to assume that the sequence\ $(C_{\nu})_{\nu \in \mathbb{N}}$ \ converges\footnote{\ $C$ \ may be the empty \ $n-$cycle.} to \ $C \in \mathcal{C}_n^{loc}(M)$. Then \ $(\sigma, C)$ \ lies in \ $F$ \ and 1) is proved. But 2) is already a consequence of the compactness of \ $p^{-1}(L)$, so \ $p : \Gamma \to S$ \ is proper.\\
To prove  that \ $p : \Gamma \to S$ \ is a topological modification of \ $S$ \ along \ $\Sigma$, it is now enough to prove that \ $p^{-1}(\Sigma)$ \ has no interior point in \ $\Gamma$. But this is obvious from the density of the graph of \ $\varphi$ \ in \ $\Gamma$. To conclude case i), notice that the projection \ $q : \Gamma \to \mathcal{C}_n^{loc}(M)$ \ is continuous and extends \ $\varphi$.\\

We shall try to explain in the comment following theorem \ref{f-mero. quotient}  why  in this case  the locally compact subset \ $\Gamma$ \ which is proper on \ $S$ \  may not be, in general,  a finite dimensional complex space.\\

To prove  ii) we first notice that the assumption that the map \ $\varphi^f$ \ is properly extendable along \ $\Sigma$ \ is equivalent to the fact that \ $\Gamma^f$, the closure of the graph of \ $\varphi^f$ \ in \ $S \times \mathcal{C}_n^f(M)$, is proper on   \ $S$. This is proved in the following lemma.

\begin{lemma}
The map \ $\varphi^f$ \ is properly extendable along \ $\Sigma$ \ if and only if the closure \ $\Gamma^f$ \  in \ $S \times \mathcal{C}_n^f(M)$ \ of the graph of \ $\varphi^f$ \ is\ $S-$proper.
\end{lemma}

\parag{Proof} Of course the properness of \ $\Gamma^f$ \ on \ $S$ \ gives a  topological modification of \ $S$ \ along \ $\Sigma$ \ with a continuous extension for \ $\varphi^f$ \ given by the projection of \ $\Gamma^f$ \ on \ $\mathcal{C}_n^f(M)$.
Conversely, if we have a proper topological modification \ $\tau : Y \to S$ \ and a continuous map \ $\psi : Y \to \mathcal{C}_n^f(M)$ \ extending \ $\varphi^f$, let  \ $\tilde{\Gamma}$ \ be the graph of \ $\psi$. Then \ $(\tau\times Id)(\tilde{\Gamma})$ \ is obviously proper on \ $S$ \ and contained in \ $\Gamma^f$, the closure of the graph of \ $\varphi^f$. But the continuity of \ $\psi$ \ implies that  \ $\tilde{\Gamma} = \Gamma^f$, as \ $\tau^{-1}(\Sigma)$ \ has no interior point in \ $Y$. $\hfill \blacksquare$

\bigskip

So to finish the proof of  ii), it is enough to endow the locally compact topological space \ $\Gamma^f$ \ with a natural structure of a weakly normal complex space such that its projection on \ $S$ \ becomes holomorphic. This is done by induction on the dimension of the "big" fibers of \ $M$ \ on \ $S$ \ in [M.00]. $\hfill \blacksquare$

\parag{Remarks}
 The key points of the proof of [M.00] are the following facts :
\begin{itemize}
\item The properness condition on \ $\Gamma^f$ \ is local along \ $\Sigma$ \ and invariant by local (proper) modification on \ $S$.
\item The quasi-properness of \ $f : M \to S$ \ is preserved by local proper modification along \ $\Sigma$ \ because of the assumption on \ $\Gamma^f$. This may not be true in presence of escape at infinity for a limit of generic fibers, so without the properness of \ $\Gamma^f$ \ on \ $S$. For the convenience of the reader, we give in the next proposition a proof of this fact as  it is a key point in the induction for proving the geometric f-flattening theorem.
\item Because of the quasi-properness of \ $f$ \ (which stays in the induction), to decrease the dimension of the biggest fibers, a local blowup of an analytic subset of \ $\Sigma$ \ is enough, as all irreducible components of all nearby  fibers meet a given compact set. Of course, with infinitely many components of dimension \ $> n$ \  in a fiber, this argument  would not work.
\item The Kuhlmann's theorem is generalized to the case of a semi-proper holomorphic map with values in an open set of a Banach space in  [M.00]. We explicit  in the appendix, for the convenience of the reader, how this generalization  is used to have a semi-proper direct image theorem with values in \ $\mathcal{C}_n^f(M)$ \ in a sense which is precised there. This formulation of the semi-proper direct image theorem with values in \ $\mathcal{C}_n^f(M)$ \ is interesting by itself. Of course an easy corollary is the "universal reparametrization theorem" of [M.00]. \\
\end{itemize}

\begin{prop}[Stability of strong quasi-properness by modification.]
Let \ $f : M \to S$ \ be a quasi-proper holomorphic map between two reduced irreducible complex spaces. Let \ $\Sigma \subset S$ \ be a closed analytic subset with no interior point  such that  the restriction of  \ $f$ \ to \ $M \setminus f^{-1}(\Sigma)$ \ is geometrically f-flat. Let \ $\varphi^f : S \setminus \Sigma \to \mathcal{C}_n^f(M)$ \ be the classifying map of the f-analytic family of \ $n-$cycles in \ $M$ associated to the fibers of  this restriction of \ $f$. Assume that \ $\Gamma^f$, the closure in \ $S \times \mathcal{C}_n^f(M)$ \ of the graph of \ $\varphi^f$, is \ $S-$proper. Let \ $\tau : \tilde{S} \to S$ \ be a proper modification of \ $S$ \ with center contained in \ $\Sigma$ \ and let \ $\tilde{M}$ \ be the strict transform of \ $M$ \ by \ $\tau$ \ and \ $\tilde{f} : \tilde{M} \to \tilde{S}$.\\
Then \ $\tilde{f}$ \ is quasi-proper and if \ $\tilde{\Sigma} : = \tau^{-1}(\Sigma)$, the closure \ $\tilde{\Gamma}^f$ \ in \ $\tilde{S} \times \mathcal{C}_n^f(\tilde{M})$ \ of the graph of the classifying map \ $\tilde{\varphi}^f : \tilde{S} \setminus \tilde{\Sigma} \to \mathcal{C}_n^f(\tilde{M})$ \ is proper on \ $\tilde{S}$.
\end{prop}

The condition  that the center of \ $\tau$ \ is contained in \ $\Sigma$ \ does not reduce the generality of the statement because f-geometric flatness is invariant by pull-back.

\parag{Proof} Recall first that the strict transform \ $\tilde{M}$ \ of \ $M$ \ by \ $\tau$ \ is the irreducible component of the fiber product \ $\tilde{S}\times_S M$ \ which dominates \ $M$. As \ $\tau$ \ is proper, the projection \ $\pi : \tilde{M} \to M$ \ is also proper and it is a proper modification along \ $f^{-1}(\Sigma)$. To show that \ $\tilde{f} : \tilde{M} \to \tilde{S}$ \ is quasi-proper, take \ $\tilde{s}_0 \in \tilde{S}$ \ and let \ $s_0 = \tau(\tilde{s}_0)$. The quasi-properness of \ $f$ \ gives an open set \ $S'$ \ in \ $S$ \ containing \ $s_0$ \ and a compact set \ $K$ \ in \ $M$ \ such that any irreducible component of \ $f^{-1}(s)$ \ for \ $s \in S'$ \ meets \ $K$. Let us show that any irreducible component of \ $\tilde{f}^{-1}(\tilde{s})$ \ for \ $\tilde{s} \in \tilde{S}' : = \tau^{-1}(S')$ \ meets the compact \ $(\tau^{-1}(f(K))\times K) \cap \tilde{M}$.\\
Let  \ $C$ \ be an irreducible component of \ $\tilde{f}^{-1}(\tilde{s})$ \ for some \ $\tilde{s} \in \tilde{S}'$, and choose a smooth point \ $(\tilde{s},x)$ \ of this fiber belonging to \ $C$. Let \ $(s_{\nu}, x_{\nu})_{\nu \in \mathbb{N}}$ \ be a sequence of points in \ $M \setminus f^{-1}(\Sigma)$ \ converging to \ $(\tilde{s},x)$ \ in \ $\tilde{M}$. Up to choose a subsequence, we may assume that the fibers \ $f^{-1}(s_{\nu})$ \ converge to a \ $n-$cycle \ $C_0$ \ in \ $\mathcal{C}_n^f(M)$, thanks to the properness of \ $\Gamma^f$ \ on \ $S$. Now we have the inclusion \ $\{\tilde{s}\}\times \vert C_0\vert \subset \tilde{f}^{-1}(\tilde{s})$, and \ $x$ \ belongs to some irreducible component \ $C_1$ \ of \ $C_0$. Then we have \ $C_1 = C$, because \ $C$ \ is the only irreducible component of \ $\tilde{f}^{-1}(\tilde{s})$ \ containing the point \ $(\tilde{s},x)$.\\
But now, we know that \ $C_1$ \ meets \ $K$. So \ $C$ \ meets \ $\tilde{K}$.\\

The proof that \ $\tilde{\Gamma}^f$ \ is proper over \ $\tilde{S}$ \ is easy : the generic fibers of \ $\tilde{f}$ \ correspond to the generic fibers of \ $f$ \  via the direct image map \ $\pi_*$. So we have a commutative diagram of continuous maps
$$ \xymatrix{\tilde{\Gamma}^f \ar[d]^{\tilde{p}_1} \ar[r]^{\pi_*} & \Gamma^f \ar[d]^{p_1} \\
\tilde{S} \ar[r]^{\tau} & S } $$
with \ $\tau$, $p_1$ \ and \ $\pi_*$ \ proper. So \ $\tilde{p}_1$ \ is proper. $\hfill \blacksquare$

\parag{Example} Let \ $f : M \to S$ \  be a quasi-proper $n-$equidimensional map between two reduced and irreducible complex spaces, where \ $n : = \dim M - \dim S$. Assume \ $S$ \ is normal so that \ $f$ \ is a geometrically f-flat map. Assume that at the generic points of  a closed irreducible curve \ $C \subset S$ \ the fibers of \ $f$ \ have two distinct  reduced and  irreducible components. Denote them by \ $X_c$ \ and \ $Y_c$ \ for \ $c \in C$. Fix a generic  point \ $c_0 \in C$ \ such that \ $X_c$ \ converges to  \ $X_{c_0}$ \ when \ $c \to c_0$.\\
Define now \ $M' : = M \setminus X_{c_0}$. Then the induced map \ $f' : M' \to S$ \ is still geometrically flat  (but not f-flat, as we shall see), because we keep the equidimensionality for \ $f'$ \ and the normality of \ $S$.  For \ $s \in S$ \ the fiber of \ $f'$ \ has only finitely many irreducible components, so the classifying map \ $\varphi : S \to \mathcal{C}_n^{loc}(M')$ \  for the fibers of the map \ $f'$ \ factors set theoretically  through the inclusion \ $\mathcal{C}_n^f(M') \hookrightarrow \mathcal{C}_n^{loc}(M')$. Let  \ $\varphi' : S \to \mathcal{C}_n^f(M')$ \ be the corresponding map. \\
Now we shall show that the closure of the graph of \ $\varphi'$ \ in \ $S \times \mathcal{C}_n^f(M')$ \ is not proper on \ $S$. Assume the contrary. Then there exists an open set \ $S' \subset S$ \ contaning \ $c_0$ \  and a compact set \ $K$ \ in \ $M'$ \ such that for any \ $c \in S', c \not= c_0$, we have \ $X_c \cap K \not= \emptyset$. Let \ $(c_{\nu})_{\nu \geq 1}$ \ be a sequence in \ $C \setminus \{c_0\}$ \ converging to \ $c_0$. Then for \ $\nu \gg 1$ \ we have \ $c_{\nu} \in S'$ \ and so \ $X_{c_{\nu}} $ \ meets \ $K$. So, up to pass to a subsequence, we may assume that we have a converging sequence \ $(x_{\nu})_{\nu \in \mathbb{N}}$ \ in \ $K$ \ to a point \ $x \in K$ \ and such that \ $x_{\nu} \in X_{c_{\nu}} $. But as the cycles \ $X_{c_{\nu}}$ \ converge to the empty $n-$cycle in \ $M'$, we obtain a contradiction.\\
Of course this shows that \ $\varphi'$ \ is not continuous. But remark that the map \ $\varphi : S \to \mathcal{C}_n^{loc}(M')$ \ is continuous (it classifies an analytic family of cycles) so its graph is closed in \ $S \times \mathcal{C}_n^{loc}(M')$ \ and its projection on \ $S$ \ is an homeomorphism (so it is proper). \\
It is easy to see that for any local  proper holomorphic modification \ $\tau : \tilde{S}' \to S'$ \ where \ $S'$ \ is an open neighbourhood of \ $c_0$, the strict transform of \ $M'$ \ cannot be geometrically f-flat on \ $\tilde{S}'$. \\
A concrete example is given as follows : let  \ $M : = \{(t,x) \in \C^2 \ / \  x^2 = t \}, S : =  \C$ \ and \ $f(t,x) = t$. Then \ $f' : M' : = M \setminus \{(1,1)\} \to S$ \ is an example of a geometrically flat map where the classifying maps
$$ \varphi : S \to \mathcal{C}_0^{loc}(M')\quad {\rm and} \quad   \varphi' : S \to \mathcal{C}_0^{f}(M')$$
have  closed graphs which are  respectively proper and not proper on \ $S$. Of course \ $\varphi'$ \ is not continuous at \ $t = 1$ \ and the map \ $f'$ \ is not quasi-proper.

\subsection{Geometric f-flattening.}

As it becomes apparent, the notion of quasi-properness  is not strong enough in presence of big fibers for an holomorphic map \ $f : M \to S$ \ between reduced irreducible complex spaces : the irreducible components of a big fibers may meet a compact set but limits of generic fibers may have an irreducible component escaping at infinity inside these big fibers.\\

These considerations lead to the following definition.

\begin{defn}\label{strongly quasi-proper}
Let \ $f : M \to S$ \ be a quasi-proper  holomorphic map between two reduced and irreducible complex spaces. We shall say that \ $f$ \ admits  {\bf local  geometric f-flattenings} if for each point in \ $S$ \ there exists an open neighbourhood \ $S'$ \ and a  proper modification \ $\tau' : \tilde{S'} \to S'$ \ such that the strict transform \ $\tilde{M}' $ \ of \ $M' = f^{-1}(S')$ \ has a geometric f-flat projection \ $\tilde{f}' : \tilde{M}' \to \tilde{S}' $. So  the following diagram commutes
$$\xymatrix{\tilde{M}' \ar[d]^{\tilde{f}'} \ar[r]^{\pi'} & M' \ar[d]^{f'} \\ \tilde{S}' \ar[r]^{\tau'} & S' } $$
and \ $\pi'$ \ is a proper modification.
In this situation we  shall say that such a quasi-proper map \ $f$ \ is {\bf strongly quasi-proper}.
\end{defn}

From the previous theorem, a quasi-proper map \ $f$ \ admits local geometric flattenings if and only if the classifying map \ $\varphi^f$ \ is properly extendable along \ $\Sigma$, and in this case, there exists a global quasi-proper f-flattening for \ $f$.\\

Coming back to the point of view of families of finite type \ $n-$cycles, this gives a good definition for f-meromorphic families :

\begin{defn}\label{f-mero.}
Let \ $S$ \ be a reduced an irreducible complex space and  let \ $\Sigma \subset S$ \ be a closed analytic subset with no  interior point  in \ $S$. Let \ $M$ \ be a reduced complex space. 
We shall say that an f-analytic  family \ $(X_s)_{s \in S \setminus \Sigma}$ \ of \ $n-$cycles of \ $M$ \ is  \ {\bf f-meromorphic along \ $\Sigma$} \ if there exists, for each point in \ $\Sigma$, an open neighbourhood \ $S'$, a proper modification \ $\tau : \tilde{S}' \to S'$ \ and a f-analytic family \ $(X_{\tilde{s}})_{\tilde{s} \in \tilde{S}'}$ \ extending the restriction of the given family to \ $S' \setminus \Sigma$.
\end{defn}

Then we have the following necessary and sufficient condition for meromorphy along \ $\Sigma$

\begin{prop}\label{f-mero}
A necessary and sufficient condition for an  f-analytic  family \ $(X_s)_{s \in S \setminus \Sigma}$ \ of \ $n-$cycles of \ $M$ \ to be f-meromorphic along \ $\Sigma$ \ is that there exists a closed analytic subset \ $\tilde{G} \subset S \times M$ \ with the following properties :
\begin{enumerate}[i)]
\item \ $\tilde{G} \cap (S \setminus \Sigma)\times M = \vert G\vert $ \ where \ $G$ \ is the graph of the given family.
\item The projection \ $ \tilde{p} : \tilde{G} \to S$ \ is strongly quasi-proper.
\end{enumerate}
\end{prop}

\parag{Proof} The condition is clearly necessary, because the image of the graph of the "extended" family by the map \ $(\tau, \pi) : \tilde{S} \times \tilde{M} \to S \times M$ \ is proper.\\
The condition is sufficient, because the existence of \ $\tilde{G}$ \ shows that the closure of the graph of classifying map of the given family is proper on \ $S$. Then we may apply the f-flattening theorem \ref{flattning 1} and use a  normal \ $\tilde{S}$ \ to conclude the f-analyticity of the extended family . $\hfill \blacksquare$\\

So the most significative part of the f-flattening theorem \ref{flattning 1} may be reformulate in the following way

\begin{thm}[f-flattening, second version]
Let \ $M$ \ and \ $S$ \ be reduced irreducible complex spaces and \ $(X_s)_{s \in S \setminus \Sigma}$ \ be an f-analytic family of cycles in \ $M$ \ parametrized by \ $S \setminus \Sigma$, where \ $\Sigma$ \ is a closed analytic subset with no interior points in \ $S$. Assume that the support of the graph of this family is the restriction to \ $(S \setminus \Sigma) \times M$ \ of an irreducible analytic subset  \ $X$ \ of \ $S \times M$ \ which is quasi-proper over \ $S$. Then the following conditions are equivalent :
\begin{enumerate}
\item The quasi-proper map \ $p : X \to S$ \ admits local f-flattenings along \ $\Sigma$ ;
\item The quasi-proper map \ $p : X \to S$ \ admits a global  f-flattening along \ $\Sigma$ ;
\item The family \ $(X_s)_{s \in S \setminus \Sigma}$ \ of \ $n-$cycles of \ $M$ \ is f-meromorphic along \ $\Sigma$.
\item The classifying map \ $\varphi^f : S \setminus \Sigma \to \mathcal{C}_n^f(M)$ \ of the family  \ $(X_s)_{s \in S \setminus \Sigma}$ \ is properly extendable along \ $\Sigma$.
\item There exists a proper (holomorphic) modification \ $\tau : \tilde{S} \to S$ \ with center contained in \ $\Sigma$ \ and a f-analytic family \ $(X_{\tilde{s}})_{\tilde{s} \in \tilde{S}}$ \ of \ $n-$cycles in \ $M$ \ extending the given f-analytic family on \ $S \setminus \Sigma$.
\end{enumerate}
\end{thm}

Now using the reparametrization theorem for a global f-flattening of the map \ $f$ \ we find a {\bf canonical modification} \ $\tau : \tilde{S} \to S$ \ where \ $\tilde{S}$ \ is simply the closure of the graphe of the map \ $\varphi^f$ \ in \ $S \times \mathcal{C}_n^f(M)$ \ endowed with a structure of a weakly normal complex space, thanks to the semi-proper direct image theorem of D. Mathieu (used in a Banach analytic setting, see [M.00] or the appendix).

\begin{defn}
 The  geometrically f-flat map \ $\tilde{f} : \tilde{M} \to \tilde{S}$ \ obtained in this way, where \ $\tilde{M}$ \ is the strict transform by \ $\tau$ \ of \ $M$, will be called the {\bf canonical f-flattening} of \ $f$.
 \end{defn}
 
 We let the reader give the  obvious universal property of the canonical f-flattening.
 
 \begin{prop}\label{big fiber}
 Let \ $f : M \to S$ \ be a strongly quasi-proper holomorphic map between reduced and irreducible complex spaces. Let \ $\tau : \tilde{S} \to S$ \ be a proper modification of \ $S$ \ such that the strict transform \ $\tilde{f} : \tilde{M} \to \tilde{S}$ \ is geometrically f-flat. Then for any \ $s \in S$ \ we have
 $$f^{-1}(s) = \bigcup_{\tau(\tilde{s}) = s} \ p\big(\tilde{f}^{-1}(\tilde{s})\big) ,$$
 where \ $p : \tilde{M} \to M$ \ is the projection.
 \end{prop}
 
 This proposition shows that in a strongly quasi-proper map the limits of generic fibers in the sense of the topology of \ $\mathcal{C}_n^f(M)$) fill up the big fibers.
 
 \parag{Proof} Let \ $\vert \tilde{G}\vert$ \ be the graph of the f-analytic family \ $(X_{\tilde{s}})_{\tilde{s} \in \tilde{S}}$ \ extending the analytic family of fibers of \ $S$ \ on \ $S \setminus \Sigma$. Then the holomorphic map \ $\tau\times id_M : \vert \tilde{G}\vert \to M$ \ is proper so for each \ $s \in S$ \ the image of \ $\tilde{f}^{-1}(\tau^{-1}(s))$ \ which is equal to  \ $ \cup_{\tau(\tilde{s})} \  p\big(\tilde{f}^{-1}(\tilde{s})\big)$ \ is a closed analytic subset in \ $f^{-1}(s)$. \\
 If \ $(s,x)\in f^{-1}(s)$ \ is not in this subset, the point \ $x \in M$ \ has an open neighbourhood \ $V$ \ which does not meet \ $X_{s'}$ \ for any \ $s' \in S\setminus \Sigma$. But then an irreducible component of \ $M$ \ containing a non empty open set in  \ $V$ \ does not meet \ $X_{s'}$ \ for any \ $s' \in S\setminus \Sigma$ \ for \ $s'$ \ near enough to \ $s$. As \ $M$ \ is irreducible, this is impossible. $\hfill \blacksquare$

 \parag{Remark} Let \ $S$ \ be a compact topological space and let \ $(X_s)_{s \in S}$ \ be a continuous family of \ $n-$cycles in a reduced complex space \ $M$. Then the projection of the set-theoretic graph \ $p : \vert G\vert \to M$ \ of the family is  proper on \ $M$. The part i) of the first version of the flattening theorem \ref{flattning 1}, gives always a proper topological  modification; so we may conclude that the union of limits of generic fibers, in the sense of the topology of \ $\mathcal{C}_n^{loc}(M)$, in a big fiber is a closed set which is a union of closed analytic \ $n-$dimensional subsets. Now the same argument as in the proof above shows that, assuming again that  \ $M$ \ is irreducible,  the big fibers are filled up by limits in the topology of \ $\mathcal{C}_n^{loc}(M)$,  of generic fibers  \\
 But notice that we have more limits for this topology than for the topology of \ $\mathcal{C}_n^f(M)$. And the assertion would not be true with this topology even if we assume \ $f$ \ quasi-proper because in this case the quasi-properness of \ $f$ \ is not enough, in general, to show that the closure of the graph of the family of generic fibers is proper on \ $S$.

 \section{Quasi-proper meromorphic quotients.}

\subsection{Quasi-proper analytic equivalence relations.}

Let now \ $X$ \ be an irreducible complex space and \ $R \subset X \times X$ \ an analytic subset which is the graph of an equivalence relation denoted  by \ $\mathcal{R}$. 

\begin{defn}
We shall say that \ $\mathcal{R}$ \ is a {\bf quasi-proper (resp. strongly quasi-proper) analytic equivalence relation} on \ $X$ \ if the following conditions hold
\begin{enumerate}[i)]
\item  \ $R$ \ has a finite number of irreducible components \ $(R_i)_{i \in I}$.
\item The projection \ $\pi_i : R_i \to X$ \ (on the first factor) is quasi-proper for each \ $ i \in I$.
\item All irreducible components of \ $R$ \ which surject on \ $X$ \ have the same dimension \ $\dim X + n$ \ (and are strongly quasi-proper on \ $X$).
\end{enumerate}
\end{defn}

In this situation the irreducible components of \ $R$ \ which does not surject on \ $X$ \ will be forgotten and we shall denote by \ $R_s$ \ their union and  by \ $R_b$ \ the union of the others components. Using the direct image theorem of Kuhlmann gives the following facts
\begin{enumerate}
\item The projection on \ $X$ \ of \ $R_s$ \ is an analytic set with empty interior.
\item The set in \ $X$ \ where the dimension of  \ $\pi^{-1}(x)$ \ is \ $> n$ \ is an analytic set in \ $X$ \ with no interior point.
\end{enumerate}

\parag{Notation}
We shall denote by \ $\Sigma$ \ the analytic subset in \ $X$ \ which is the union of the non normal points in \ $X$, the image of \ $R_s$ \ in \ $X$ \ and the set in \ $X$ \ where the dimension of  \ $\pi^{-1}(x)$ \ is \ $> n$. It has  no interior point in \ $X$, and the holomorphic map 
$$ \pi : R \setminus \pi^{-1}(\Sigma) \to X \setminus \Sigma $$
is a quasi-proper surjective \ $n-$equidimensional holomorphic map on a normal complex space. So we have an f-analytic family of \ $n-$cycles in \ $X$ \ associated to its fibers and so a classifying map
$$ \varphi : X \setminus \Sigma \to \mathcal{C}_n^f(X) $$
for this family.

\begin{defn}\label{mero. quotient}
In the situation describe above we shall say that a meromorphic map
$$q :  X \dashrightarrow Q $$
is a { \bf quasi-proper meromorphic quotient} (resp. {\bf strongly quasi-proper meromorphic quotient}) for the equivalence relation \ $\mathcal{R}$ \ when there exists a proper modification \ $\tau : \tilde{X} \to X$ \ with center contained in \ $\Sigma$ \ and a holomorphic map \ $\tilde{q} : \tilde{X} \to Q$ \ inducing the meromorphic map \ $q$ \  (note that this implies that \ $q$ \ is holomorphic on \ $X \setminus \Sigma$) such that the following conditions are satisfied :
\begin{enumerate}[i)]
\item For \ $x, y \in X \setminus  \Sigma$ \ we have \ $x \mathcal{R} y$ \ if and only if \ $q(x) = q(y)$.
\item The map \ $\tilde{q}$ \ is quasi-proper (resp. strongly quasi-proper) and surjective.
\end{enumerate}
When \ $R$ \ is strongly quasi-proper on \ $X$ \  we shall say that the meromorphic map \ $q : X \dashrightarrow Q$ \ is the  { \bf universal meromorphic quotient} when the corresponding holomorphic map \ $\tilde{q} : \tilde{X} \to Q$ \ is the universal reparametrization of the family of fibers of the canonical f-flattening of the strongly quasi-proper map \ $p_1 : R_b \to X$.
\end{defn}

\parag{Remarks}
\begin{enumerate}
\item It is easy to see that  if  any strongly  quasi-proper quotient exists, then the universal meromorphic quotient also exists and they are  meromorphically equivalent.
\item  For a given quasi-meromorphic quotient \ $q : X \dashrightarrow Q$ \ the fact that after some proper modification the map \ $\tilde{q}$ \ becomes strongly quasi-proper is independant of the choice of the modification \ $\tau : \tilde{X} \to X$ \ on which \ $q$ \ extends holomorphically because this property only depends on generic fibers of \ $q : X \setminus \Sigma \to \mathcal{C}_n^f(M)$.\\
\end{enumerate}

\subsection{Existence of a meromorphic quotient.}

Our main result is the following theorem.

\begin{thm}\label{f-mero. quotient}
Let \ $X$ \ be a reduced irreductible complex space and \ $\mathcal{R}$ \ be a strongly quasi-proper analytic equivalence relation, so  satisfying the following condition :
\begin{itemize}
\item Each irreducible component of the graph \ $R$ \ of \ $\mathcal{R}$ \ which surjects on \ $X$ \ is strongly quasi-proper on \ $X$ (or admits a geometric  f-flattening on \ $X$). $\hfill (@) $
\end{itemize}
Then \ $\mathcal{R}$ \ admits a universal strongly quasi-proper meromorphic quotient. The set \ $Q$ \ is the weak normalisation of the image in \ $\mathcal{C}_n^f(X)$ \ of the closure in \ $X \times \mathcal{C}_n^f(X)$ \ of the graph of the map \ $\varphi^f : X \setminus \Sigma \to \mathcal{C}_n^f(X)$ \ classifying the f-analytic family whose graph is \ $R_b \cap ((X \setminus \Sigma)\times X)$.
\end{thm}

\parag{Comment} With our method we cannot prove the existence of a quasi-proper meromorphic quotient   
in the case of a quasi-proper equivalence relation \ $\mathcal{R}$ \ which does not satisfy the condition \ $(@)$. This comes from the fact that the generalization given in [M.00] of Kuhlmann's theorem may be used only in the case of the classifying map of a f-analytic family;  in the case of a semi-proper  "holomorphic" map of a finite dimensional complex space with values in \ $\mathcal{C}_n^f(X)$ \ we can use banach analytic sets to determine (locally) finite type cycles (see the appendix); but this does not seem possible for the analoguous case with values in \ $\mathcal{C}_n^{loc}(X)$ : for the topology of \ $\mathcal{C}_n^{loc}(X)$ \ cycles near a given cycle cannot be determined by a finite number of (adapted)  scales.

\parag{Proof} Consider the  modification \ $\tau : \tilde{X} \to X$ \ given by the f-flattening theorem applied to the quasi-proper map \ $R_b \to X$ \ which admits a f-geometric flattening by assumption. Now the strict transform \ $\tilde{R}$ \ of \ $R_b$ \ via \ $\tau$ \ has a geometrically f-flat projection on \ $\tilde{X}$. So we get a classifying map
$$ \tilde{\varphi}^f : \tilde{X} \to \mathcal{C}_n^f(X).$$
Now the point is to prove that this map is semi-proper.\\
Let \ $C_0 \in \mathcal{C}_n^f(X)$; by definition of the topology of this space, for each  open set \ $W \subset\subset  X$ \ such that each irreducible components of \ $C_0$ \ meets \ $W$, the set \ $\Omega(W)$ \ of all \ $C$ \ with this property is an open set in \ $\mathcal{C}_n^f(X)$. Let us show that we have \ $\tilde{\varphi}(\tilde{X}) \cap \Omega(W) = \tilde{\varphi}(\tilde{K}) \cap \Omega(W)$, where \ $\tilde{K} = (\tau^{-1}(f(K))\times K) \cap \tilde{X}$ \ with \ $K : = \bar W$.\\
Let \ $C = \tilde{\varphi}^f(y)$ \ be in \ $\Omega(W)$ \ for some \ $y \in \tilde{X}$. Then \ $\tau(y)$ \ is the limit of a sequence of points \ $(x_{\nu})_{\nu\in \mathbb{N}}$ \ in \ $X \setminus \Sigma$. So \ $\tilde{\varphi}^f(x_{\nu})$ \ converges to \ $C$. For \ $\nu$ \ large enough  \ $\tilde{\varphi}^f(x_{\nu}) $ \ lies in \ $\Omega(W)$. This means that the equivalence class of \ $x_{\nu}$ \ which is \ $\vert \varphi^f(x_{\nu})\vert$ \ meets \ $W$. Then, for \ $\nu$ \ large enough,  there exists a point \ $x'_{\nu}\in W$ \ such \ $\tilde{\varphi}^f(x'_{\nu}) = \tilde{\varphi}^f(x_{\nu})$. Consider now the sequence \ $(x'_{\nu})_{\nu \geq \nu_0}$ \ as a sequence in \ $\tau^{-1}(\bar W)$. We may, up to pass to a subsequence, assume that it converges to a point \ $z \in \tau^{-1}(\bar W)$. By continuity of \ $\tilde{\varphi}^f$ \ we shall have \ $\tilde{\varphi}^f(z) = \varphi^f(y)$ \ and by construction, \ $z$ \ is in \ $ \tilde{K}$. \\
So the semi-properness of \ $\tilde{\varphi}^f$ \ is proved, and the "reparametrization theorem" of [M.00] (see also [B.08] or the Appendix) may be applied to \ $\tilde{\varphi}^f$. Then it gives a weakly normal  f-meromorphic quotient for \ $\mathcal{R}$ \ in the sense of the definition \ref{mero. quotient}, because our proof applies now to the map \ $\tilde{q} : \tilde{X} \to Q$ \ where \ $Q$ \ is the weak normalisation of the image of \ $\tilde{\varphi}^f$; it gives that \ $\tilde{q}$ \  is in fact geometrically f-flat. $\hfill \blacksquare$

\subsection{Extension to meromorphic equivalence relations.}

Let me give two definitions in order to formulate a simple corollary of the meromorphic quotient theorem.

\begin{defn}\label{mero. equiv.}
Let \ $X$ \ be a reduced and  irreducible complex space. A {\bf quasi-proper meromorphic equivalence  relation} \ $\mathcal{R}^m$ \  on \ $X$ \ will be a closed  analytic subset \ $R^m \subset X \times X$ \ with finitely many irreducible components, such that there exists a closed analytic subset \ $Y \subset X$ \ with no interior point in \ $X$ \ with the following conditions :
\begin{enumerate}[i)]
\item \ $R : = R^m \cap\big[ (X\setminus Y)\times(X \setminus Y)\big] $ \ is the graph of an equivalence relation \ $\mathcal{R}$ \ on \ $X \setminus Y$.
\item We have \ $\bar R = R^m$ \ in \ $X \times X$.
\item Each irreducible component of \ $R^m$ \ is quasi-proper on \ $X$ \ via the first projection.
\end{enumerate}
We shall say that \ $\mathcal{R}^m$ \ is  {\bf strongly quasi-proper}  when there exists an integer \ $n $ \ such that  the first projection of each irreducible component of \ $R$ \ which surjects on  \ $X$ \ is  strongly quasi-proper with generic fiber of pure dimension \ $n$.
\end{defn}

\parag{Remark} If,  for a meromorphic equivalence relation on \ $X$,  the subset  \ $R^m$ \ is proper on \ $X$ \ via the first projection,  we shall say that the meromorphic equivalence relation \ $\mathcal{R}^m$ \ is proper. Of course, in this case, the geometric flattening theorem for compact cycles (see [B.78]) implies that \ $\mathcal{R}^m$ \ is strongly quasi-proper as soon as \ $R$ \ is generically equidimensional on \ $X$.

\bigskip

\begin{defn}\label{mero.quot.}
In the situation of the previous definition, we shall say that \ $\mathcal{R}^m$ \ has a {\bf quasi-proper (resp. strongly quasi-proper, resp. proper) meromorphic quotient} when there exists a  meromorphic surjective map \ $ q : X \dashrightarrow Q$ \ whose graph is quasi-proper (resp. strongly quasi-proper, resp. proper) on \ $Q$ \ such that there exists a dense  open set \ $\Omega$ \ of  \ $X \setminus Y$, on which  \ $q$ \ is holomorphic and such that  for \ $x,x'$ \ in \ $\Omega$ \ we have \ $x\mathcal{R} x' $ \ if and only if \ $q(x) = q(x')$.
\end{defn}

\begin{thm}
Let \ $X$ \ be a reduced irreductible complex space and \ $\mathcal{R}^m$ \ be a  strongly quasi-proper  (resp. proper) meromorphic equivalence relation on \ $X$. Then \ $\mathcal{R}^m$ \ admits a strongly quasi-proper (resp. proper) meromorphic quotient. 
\end{thm}

\parag{Proof} Define \ $R_b$ \ as the union of the components of \ $R^m$ \ which dominates \ $X$ \ and \ $R_s$ \ as the union of the other components. Let \ $Z$ \ be the projection in \ $X$ \ of \ $X_s$. It is a closed analytic subset with no interior points. We choose \ $\Sigma \subset X$ \ such that it contains \ $Z$. non normal points of \ $X$, the set of points \ $x \in X$ \ such that the fiber at \ $x$ \ of the projection \ $R_b \to X$ \ is \ $> n$ \ and also \ $Y$. Then we have a classifying map
$$ \varphi^f : X \setminus \Sigma \to \mathcal{C}_n^f(X) $$
and the closure \ $\Gamma^f$ \ of its graph in \ $X \times \mathcal{C}_n^f(X)$ \ is proper on \ $X$, by the strongly proper assumption of the projection \ $R_b \to X$. But now for \ $x \not\in \Sigma$ \ the set \ $q(x)$ \ may not be equal to its equivalence class in \ $X \setminus Y$, which is equal to \ $q(x) \cap (X \setminus Y)$. But for \ $x,x' \not\in \Sigma$ \ the equality \ $q(x) = q(x')$ \ either implies \ $x \mathcal{R} y$ \ when the two fibers of \ $R$ \ over \ $x, x'$ \ are equal when cut with \ $X \setminus Y$, or at least one has an irreducible component contained in \ $Y$. So it is enough to delete in \ $X \setminus \Sigma$ \ the analytic subset of \ $x$ \ such this happens ; this is a closed analytic subset thanks to the lemma \ref{facile};  then we  define the dense open set  \ $\Omega \subset X \setminus \Sigma$ \ as its complement. $\hfill \blacksquare$

\begin{lemma}\label{facile}
Let \ $(X_s)_{s \in S}$ \ be a f-analytic family of cycles in \ $M$ \ and \ $Y \subset M$ \ a closed analytic subset. The set of \ $s \in S$ \ such \ $X_s$ \ has an irreducible component contained in \ $Y$ \ is a closed analytic subset in \ $S$.
\end{lemma}

\parag{Proof} This set is clearly closed because a limit of cycles contained  in \ $Y$ \ is contained in \ $Y$. The problem is then local on \ $S$. But, thanks to the quasi-properness on \ $S$ \ of the graph of the family, it is enough to consider finitely many scales ;  this reduces to the analoguous lemma for multiform graphs, so when \ $S = H(\bar U, \Sym^k(B))$ \ and \ $Y \subset U\times B = M$. This case is elementary. $\hfill \blacksquare$\\

\parag{Remark} In the case of a proper meromorphic equivalence relation which has pure dimensional generic fibers we obtain that there always exists an universal  proper meromorphic quotient. Compare with [C.60].

\section{Stein factorization.}

The aim of this section is to apply the previous methods to build a Stein's factorization for strongly quasi-proper map. We shall begin by the geometric f-flat case.

\subsection{The geometric f-flat case.}

This paragraph is devoted to prove the following result :

\begin{thm}\label{Stein f-flat}
Let \ $f : M \to S$ \ be a geometrically f-flat holomorphic surjective map between two reduced  complex spaces. Assume that \ $M$ \ is normal and  denote  \ $n : = \dim M - \dim S$.  Then there exists a geometrically f-flat holomorphic map \ $g : M \to T$ \ on a weakly normal complex space \ $T$ \ and a proper and finite surjective map \ $p : T \to S $ \ such that \ $f = p\circ g$ \ and such that the generic fiber of \ $g$ \ is irreducible.
\end{thm}

We shall begin by some lemmas.

\begin{lemma}\label{irred.}
Let \ $M$ \ be a reduced complex space. The subset \ $Irr_n(M)$ \ of cycles in \ $\mathcal{C}_n^f(M)$ \ which are reduced (all multiplicities are \ $1$) and irreducible is an open set in \ $\mathcal{C}_n^f(M)$.
\end{lemma}

Notice that this lemma is false in \ $\mathcal{C}_n^{loc}(M)$.

\parag{Proof} Let \ $C_0 \in \mathcal{C}_n^f(M)$ \ be reduced and irreducible. Choose a smooth point \ $x_0$ \  in \ $C_0$ \ and an \ $n-$scale \ $E : = (U,B,j)$ \ on \ $M$ \ adapted  to \ $C_0$ \ whose center  \ $C(E) : =  j^{-1}(U \times  B)$ \ contains \ $x_0$ \ and such that \ $deg_E(C_0) = 1$. Now consider the open set \ $\Omega_1(E) \cap \Omega(C(E))$ \  in \ $\mathcal{C}_n^f(M)$, where \ $\Omega_1(E) $ \ is the set of cycles for which \ $E$ \ is adapted and have degree 1 in \ $E$, and where \ $\Omega(W)$ \ is  the set of cycles such that each irreducible component meets \ $W$. This  open set   contains \ $C_0$. Now each cycle \ $C$ \ in this open set is clearly reduced and irreducible. $\hfill \blacksquare$

 \parag{Notation} For a reduced complex space \ $M$ \ denote by \ $red\mathcal{C}_n^f(M) \subset \mathcal{C}_n^f(M)$ \ the subset of reduced cycles of \ $M$.\\

\begin{cor}\label{reduced}
Let \ $M$ \ be a reduced complex space. The subset \ $red\mathcal{C}_n^f(M)$ \ of reduced cycles in \ $\mathcal{C}_n^f(M)$ \ is open.
\end{cor}

\parag{Proof}  Let \ $C_0 \in red\mathcal{C}_n^f(M)$. Denote by \ $J$ \ the finite set of irreducible components of \ $C_0$. Choose for each irreducible component \ $C_0^j, j \in J$, of \ $C_0$ \ a smooth point \ $x^j$ \  in \ $C_0$ \ which belongs to \ $C_0^j$. Choose then, for each  \ $j \in J$, an \ $n-$scale \ $E_j$ \ adapted to \ $C_0$ \ such that \ $x^j \in C(E_j)$ \ and \ $deg_{E_j}(C_0) = 1$. Then any cycle  \ $C$ \ in the open set 
$$ \bigcap_{j \in J} \ \Omega_1(E_j) \bigcap \Omega(\cup_{j \in J} C(E_j)) $$
 is reduced. $\hfill \blacksquare$\\

\begin{lemma}\label{non reduced}
Let \ $M$ \ be a reduced complex space. The subset \ $\Sigma$ \ of cycles in \ $\mathcal{C}_n^f(M)$ \ which have a non reduced irreducible component (so at least one multiplicity is \ $\geq 2$) is a closed analytic subset in the following sense : for any f-analytic family \ $(X_s)_{s \in S}$ \ parametrized by a reduced complex space \ $S$, the pull-back \ $\varphi^{-1}(\Sigma)$ \ is a closed analytic  subset of \ $S$, where \ $\varphi$ \ is the classifying map \ $\varphi : S \to \mathcal{C}_n^f(M)$ \ of the f-analytic family.
\end{lemma}

Notice again that this lemma is false in \ $\mathcal{C}_n^{loc}(M)$.

\parag{Proof} As the complement of \ $\Sigma$ \  is open, thanks to corollary \ref{reduced},
\ $\Sigma$ \ is closed.\\
We shall now give "local holomorphic equations" for \ $\Sigma$ \ in \ $\mathcal{C}_n^f(M)$. Let \ $C_0 \in \Sigma$. Choose for each irreducible component \ $C_0^j, j \in J$, of \ $\vert C_0\vert$ \ a smooth point \ $x_j$ \ of \ $\vert C_0\vert$ in \ $C_0^j$. Let \ $E_j : = (U_j,B_j, i_j)$ \ be a \ $n-$scale adapted to \ $C_0$ \ and whose center \\ $C(E_j) : = i_j^{-1}(U_j\times B_j)$ \ contains \ $x_j$. Put \ $W : = \cup_{j \in J} \ C(E_j)$. Let \ $k_j : = deg_{E_j}(C_0)$ \ and define the open set \ $\mathcal{V}$ \ in \ $\mathcal{C}_n^f(M)$ \ as
$$ \mathcal{V} : = \big(\bigcap_{j\in J} \ \Omega_{k_j}(E_j)\big) \bigcap \Omega(W).$$
Now consider the map 
$$ \alpha : \mathcal{V} \to \prod_{j \in J} \ H(\bar U_j, \Sym^{k_j}(B_j)) $$
which associates to \ $C \in \mathcal{V}$ \ the collection of multiform graphs that \ $C$ \ defines in the adapted scales \ $(E_j)_{j \in J}$. Then \ $\alpha$ \  is continuous and injective. Continuity is obvious form the definition of the topology of \ $\mathcal{C}_n^f(M)$. Injectivity comes from the fact that for \ $C \in \mathcal{V}$ \ any irreducible component of \ $C$ \ meets \ $W$. So, if \ $\alpha(C) = \alpha(C')$ \ for \ $C, C' \in \mathcal{V}$ \ we have \ $C \cap W = C'\cap W$ \ and hence \ $C = C'$. Define \ $J' : = \{ j \in J \ / \ k_j\geq 2\}$, and  for each \ $j \in J'$ \ let \ $Z_j \subset H(\bar U_j, \Sym^{k_j}(B_j))$ \ be the subset of \ $X \in H(\bar U_j, \Sym^{k_j}(B_j))$ \ such that at least one irreducible component of \ $X \cap (U_j\times B_j)$ \ is not reduced. Then \ $Z_j$ \ is a closed banach analytic subset of \ $H(\bar U_j, \Sym^{k_j}(B_j))$ \ because it is the zero fiber of the holomorphic map given by the discriminant
$$ \Delta : H(\bar U_j, \Sym^{k_j}(B_j))  \to H(\bar U_j, S^{k_j.(k_j-1)}(\C^p)). $$
Now it is clear that the subset \ $\Sigma \cap \mathcal{V}$ \ is the pull-back by \ $\alpha$ \ of the closed banach analytic subset 
 $$\mathcal{Z} : = \bigcup_{j \in J'} \ \big(Z_j\times\prod_{j \in J\setminus\{j\}} H(\bar U_j, \Sym^k(B_j))\big) .$$
 If we have \ $\varphi : S' \to \mathcal{V}$ \ which is the classifying map of an f-analytic family of $n-$cycles in \ $M$ \ parametrized by a reduced complex space \ $S'$, the map \ $\alpha\circ \varphi$ \ is holomorphic and so the set \ $\varphi^{-1}(\Sigma) = (\varphi\circ\alpha)^{-1}(\mathcal{Z})$ \ is a closed analytic subset in \ $S'$. $\hfill \blacksquare$
 
 \begin{lemma}\label{universal irred}
 Let \ $M$ \ be a reduced complex space and \ $n$ \ an integer. Define the set \ $\Phi : = \{ (X,Y) \in \mathcal{C}_n^f(M)\times \mathcal{C}_n^f(M) \ / \ X \leq Y \}$. Then \ $\Phi$ \ is closed and its second projection on \ $\mathcal{C}_n^f(M)$ \ is proper finite and surjective.
 \end{lemma}
 
 \parag{Proof} The closedness of \ $\Phi$ \ is a consequence of the next lemma. The fact that the second projection is surjective with finite fibers is obvious. So it is enough to prove that it is a closed map. But this is an easy consequence of the description of compact sets in \ $\mathcal{C}_n^f(M)$ \ combined with the next lemma. $\hfill \blacksquare$

  \begin{lemma}\label{limite inegalite}
  Let \ $(X_{\nu})_{\nu\in \mathbb{N}}$ \ and \ $(Y_{\nu})_{\nu\in \mathbb{N}}$ \ be two converging sequences in \ $\mathcal{C}_n^f(M)$. If we have \ $X_{\nu} \leq Y_{\nu}$ \ for each \ $\nu \gg 1$ \ then the limits \ $X$ \ and \ $Y$ \ satisfy \ $X \leq Y$.
  \end{lemma}
  
  \parag{Proof} The inclusion \ $\vert X\vert \subset \vert Y\vert $ \ is clear. Let \ $C$ \ be an irreducible component of \ $Y$. Let \ $y$ \ be a smooth point of \ $\vert Y\vert$ \ in \ $C$, and let \ $E$ \ be an \ $n-$scale on \ $M$ \ whose center contains \ $y$ \ and such that \ $deg_E(C) = 1$. Let \ $k : = deg_E(Y)$. For each \ $\nu \gg 1$ \ the scale \ $E$ \ is adapted to \ $X_{\nu}$ \ and \ $Y_{\nu}$ \ and we have \ $deg_E(X_{\nu}) \leq k = deg_E(Y_{\nu})$. So we have \ $deg_E(X) \leq k$. Then the irreducible component \ $C$ \ of \ $Y$ \ has multiplicity less than \ $k$ \ in \ $X$. $\hfill \blacksquare$\\
  
  Observe that this lemma is also true for converging sequences in \ $\mathcal{C}_n^{loc}(M)$.\\

 \begin{lemma}\label{obvious 1}
 Let \ $f : M \to S$ \ be a geometrically f-flat map and let \ $N \subset M$ \ be a closed analytic subset such \ $N$ \ has empty interior in each fiber of \ $f$. Then the restriction  \ $f' : M \setminus N \to S$ \ is again a geometrically f-flat map.
 \end{lemma}
 
 \parag{Proof} If \ $(X_s)_{s \in S}$ \ is the f-analytic family of fibers of \ $f$, then \ $(X_s\setminus N)_{s \in S}$ \ is an analytic family of cycles in \ $M \setminus N$. The only point to see to conclude the proof is to show that \ $f'$ \ is quasi-proper. For any compact  \ $K$ \  in \ $S$ \ there exists a relatively compact open set \ $W_K$ \ in \ $M$ \ such that any irreducible component of any \ $X_s$ \ for \ $s \in K$ \ meets \ $W$. Let \ $L : = \bar W\cap N$, and choose a basis of  a compact neighbourhood \ $(\Lambda_{\nu})_{\nu \in \mathbb{N}}$ \ of \ $L$. If, for infinitely  many \ $\nu$,  there exists an irreducible component \ $C_{\nu}$ \ of a fiber \ $f^{-1}(s_{\nu})$ \ with \ $s_{\nu} \in K$ \ such that \ $C_{\nu}\setminus N$ \ does not meet \ $\Lambda_{\nu}$, up to pass to a subsequence we may assume that \ $s_{\nu}$ \ converges to a point \ $s \in K$ \ and  that \ $C_{\nu}$ \ converges in \ $\mathcal{C}_n^f(M)$ \ to a non empty cycle \ $C$ \ such that \ $\vert C\vert  \leq f^{-1}(s) \cap N$. This contradicts our assumption that \ $N$ \ has empty interior in each fiber of \ $f$. $\hfill \blacksquare$
 
   \parag{Proof of theorem \ref{Stein f-flat}}Define \ $\Sigma \subset S$ \ as the union of non normal points in \ $S$ \ and of  the set of points in \ $S$ \ where the fiber \ $f^{-1}(s)$ \ is not reduced (it is a closed analytic subset of \ $S$ \ thanks to lemma \ref{non reduced}). Let \ $N$ \ be the subset of points \ $x$ \  in \ $M$ \ which have multiplicity \ $\geq 2$ \ in the cycle  \ $X_{f(x)}$. From  ch.IV of [BOOK], we know that \ $N$ \ is a closed analytic subset in \ $M$. Let \ $M' : = M \setminus (N \cup f^{-1}(\Sigma))$. Then each point in \ $M'$ \ is a smooth point of the reduced cycle \ $X_{f(x)} \cap M'$. Denote by \ $C_x$ \ for \ $x \in M'$ \ the irreducible component of \ $f^{-1}(x)$ \ containing \ $x$. It is unique, by definition of \ $M'$, and we have an f-analytic family of cycles in \ $M$,  $(C_x)_{x \in M'}$, parametrized by \ $M'$  : analyticity is clear from the criterion of analytic extension for analytic  families of cycles because \ $N $ \ has empty interior in each cycle in any \ $X_s$ \ with \ $s \not\in \Sigma$ \ (see  [BOOK]) ;  quasi-properness on \ $M'$ \ comes from the quasi-properness of \ $f'$ \ proved in lemma \ref{obvious 1}. Let  
  $$ g_0 : M' \to  \mathcal{C}^f_n(M)$$
 be  the classifying map  of this family. Then the map \ $f'\times g_0$ \ takes values in the set
 $$ \hat{S} : = \{ (s, C) \in S\times \mathcal{C}_n^f(M) \ / \  C \leq X_s \} $$
 and more precisely in \ $p^{-1}_1(S \setminus \Sigma)$ \ where \ $p_1 : \hat{S} \to \mathcal{C}_n^f(M)$ \ is the first projection. Let \ $\Gamma$ \ be the image  \ $(f', g_0)(M')$ \ and \ $\bar \Gamma$ \  the closure of the image of \ $M'$ \ in \ $\hat{S}$. Note that we know, thanks to lemma \ref{universal irred},  that \ $\hat{S}$ \ is closed in \ $S \times \mathcal{C}_n^f(M)$. But as \ $f'$ \ is quasi-proper, the map \ $(f', g_0)$ \ is semi-proper, and using D. Mathieu's semi-proper direct image [M.00] (see also the appendix) we obtain that \ $\Gamma$ \ has a natural structure of weakly normal complex space ; then \ $p_1 : \Gamma \to S \setminus \Sigma$ \ is a branched covering. But using again the lemma \ref{universal irred} we obtain that \ $p_1 : \bar\Gamma \to S$ \ is also a branched covering. And now [G.R.58] (or [D.90])  gives a natural structure of weaky normal complex space on \ $\bar\Gamma$.\\
The holomorphic map \ $g' : M' \to \Gamma \subset \bar\Gamma$ \ induced by \ $(f', g_0)$ \ is locally bounded along the analytic subset \ $f^{-1}(\Sigma)\cup N$ \ in the sense that, locally in \ $S$ \ we can embed \ $\bar\Gamma$ \ in \ $S \times B$ \ where \ $B$ \ is a polydisc in an affine space, so we have, by normality of \ $M$, a holomorphic  extension \ $g : M \to \bar\Gamma$.\\
 It is now an exercice to check that \ $g$ \ is geometrically f-flat with  irreducible generic fibers,  and satifies \ $f = p_1\circ g$. So we conclude by putting \ $T : = \bar\Gamma$ \ and \ $p : = p_1$. $\hfill \blacksquare$

 \subsection{The strongly quasi-proper case.}

  Now we shall give the Stein's factorization theorem for a strongly quasi-proper surjective holomorphic map.\\
  
  The main  problem comes from the  hypothesis of normality for \ $M$ \ in our theorem \ref{Stein f-flat}. Notice first that this assumption was essential to extend to \ $M$ \ the map \ $g$ \ defined on \ $M'$.\\
  
  The problem comes from the following fact :
  \begin{itemize}
  \item  Let \ $\nu : \tilde{M} \to M$ \ be the normalisation of a reduced  irreductible complex space \ $M$. Let \ $Z \subset M$ \ be an irreducible subset. Then \ $\nu^{-1}(Z)$ \ may have infinitely many irreducible components in \ $\tilde{M}$. This is shown by the next example.
  \end{itemize}
  
  \parag{Example} Let \ $\tilde{M} : = \C^2$ \ and consider on \ $\tilde{M}$ \ the equivalence relation which identifies \ $(n,z)$ \ and \ $(-n,z)$ \ for any \ $n \in \mathbb{N}^*$ \ and any \ $z \in \C$. The quotient is a reduced and irreducible complex space of dimension 2 with normal crossing singularities along the curves \ $\Delta_n : = q(D_n)$ \ where \ $q : \tilde{M} \to M$ \ is the quotient map and \ $D_n : = \{(z_1, z_2) \in  \C^2\ / \ z_1 = n \}$ \ for \ $n \in \mathbb{N}^*$.\\
  Now let \ $Z_0 : = \{(z_1, z_2) \in \C^2 \ / \ z_1 = z_2\}$. Then \ $Z : = q(Z_0)$ \ is an irreducible curve in \ $M$ \ and \ $q^{-1}(Z) = Z \cup\big(\cup_{n \in \mathbb{N}*} \ \{(-n, n)\}\big) $. So the pull-back of the irreducible curve \ $Z$ \ in \ $M$ \ has infinitely many irreducible components in \ $\tilde{M}$. Of course  \ $q$ \ is the normalisation map for \ $M$.\\
  
  The situation is not too bad, thanks to the following simple observation.
  
  \begin{lemma}\label{finitude grosses}
  Let \ $\pi : M \to N$ \ be a proper finite and surjective holomorphic map between irreducible complex spaces. Let \ $Z \subset N$ \ be an irreducible analytic subset of dimension \ $d$. Then for any compact set \ $K \subset N$ \ meeting \ $Z$, the compact set \ $\pi^{-1}(K)$ \ meets any irreducible component of \ $\pi^{-1}(Z)$ \ which is of dimension \ $d$.
  \end{lemma}
  
  \parag{Proof} It is enough to remark that any irreducible component \ $C$ \ of \ $\pi^{-1}(Z)$ \ which is of dimension \ $d$ \ satisfies \ $\pi(C) = Z$, because the irreduciblity of \ $Z$ \ and the fact that \ $\pi(C)$ \ is an analytic subset of dimension \ $d$. $\hfill \blacksquare$\\
  
  In such a situation the union of irreducible components of \ $\pi^{-1}(Z)$ \ of dimension \ $d$ \ is a \ $d-$cycle of finite type in \ $M$. So we have in this way a map
  $$ \pi^* : \mathcal{C}_n^f(N) \to \mathcal{C}_n^f(M) .$$
  The next result gives  sufficient conditions on an f-analytic family of \ $n-$cycles in \ $N$ \ in order that the corresponding family of \ $n-$cycles in \ $M$ \ is again f-analytic.
    
  \begin{prop}\label{pull-back}
  Let \ $\pi : M \to N$ \ be a proper finite and surjective map between irreducible complex spaces. Let \ $R \subset N$ \ be the ramification set of \ $\pi$, so the minimal closed analytic subset such that \ $N \setminus R$ \ is normal and \ $\pi$ \ induced a k-sheeted (unbranched)  covering \ $M \setminus \pi^{-1}(R) \to N \setminus R$. Let \ $S$ \ be a normal complex space and \ $(X_s)_{s \in S}$ \ be a f-analytic family of \ $n-$cycles in \ $N$. Assume that the generic cycle is reduced and that  the ramification set \ $R \subset N$ \ has empty interior in the generic cycle \ $X_s$. Then there exists a unique  f-analytic family \ $(Y_s)_{s\in S}$ \ of \ $n-$cycles in \ $M$ \ such that we have \ $\pi_*(Y_s) = k.X_s$ \ for each \ $s \in S$.
  \end{prop}
  
  Observe  that, as \ $Y_s$ \ has pure dimension \ $n$, it cannot contain any irreducible components of dimension \ $< n$ \ in \ $\pi^{-1}(X_s)$.
  
  \parag{Proof} Let \ $\sigma : = \dim S$. Let \ $G \subset S \times M$ \ be the graph of the f-analytic family. Define \ $\tilde{G} \subset S \times \tilde{M}$ \ as the union of the irreducible components of \ $(id_S\times \nu)^{-1}(G)$ \ which have dimension \ $\sigma + n$. The main point is to show that \ $\tilde{G}$ \ is quasi-proper on \ $S$. As \ $G$ \ has only finitely many irreducible components and they all are of dimension \ $\sigma + n$, it is enough to prove this in the case  where \ $G$ \ is irreducible\footnote{Note that quasi-properness of \ $G$ \ implies quasi-properness of its irreducible components because we are in the strongly quasi-proper case. See [B. Mg. 10].}. Now  apply the previous lemma to the map \ $\pi : = id_S \times \nu $ \ and \ $G$. It follows that \ $\tilde{G}$ \ has only  finitely many irreducible components. Let \ $\Gamma$ \ be such a component. We want to prove that \ $\Gamma \to S$ \ is quasi-proper. But if \ $K$ \ is a compact in \ $S$, there exists a compact \ $L$ \ in \ $M$ \ such that any irreducible component of any cycle \ $X_s$ \ with \ $s \in K$ \ meets \ $L$. Then, again by the lemma, any irreducible component of a fiber\footnote{Note that \ $\Gamma$ \ has only pure dimensional fibers over \ $S$.} of \ $\Gamma$ \ over a point in \ $K$ \ has to meet the compact set \ $\Gamma \cap (K \times \nu^{-1}(L))$. To conclude, it is enough to recall that \ $\tilde{G}$ \ is the graph of an unique f-analytic family of \ $n-$cycles in \ $\tilde{M}$ \ because \ $S$ \ is normal. The fact that \ $\pi_*(Y_s) = k.X_s$ \ for generic \ $s$ \ is then obvious. The analyticity of the direct image family implies equality for any \ $s \in S$. $\hfill \blacksquare$\\
  
  The following corollary is an easy consequence of the proposition applied to the normalization map of an irreducible complex space.
  
  \begin{cor}\label{utile}
  Let \ $f : M \to S$ \ be a quasi-proper  surjective and equidimensional holomorphic map of an irreducible complex space \ $M$ \ on a normal complex space \ $S$. Then there exists a unique geometrically f-flat map \ $\tilde{f} : \tilde{M} \to S$ \ such that \ $\tilde{f} = f \circ \nu $, where \ $\nu : \tilde{M} \to M$ \ is the normalization map of \ $M$.
  \end{cor}

 Notice that \ $f$ \ is  geometrically f-flat under the hypothesis of the corollary. Observe  also that the direct image via \ $\nu$ \ of the fibers of \ $\tilde{f}$ \ are the fibers of \ $f$ \ as \ $\nu$ \ has generic degree \ $1$.

  \begin{thm}\label{Stein red.}
  Let \ $f : M \to S$ \ be a strongly quasi-proper surjective holomorphic map between two reduced complex spaces.  Let \ $n : = \dim M - \dim S$ \ and assume \ $M$ \ normal. Then there exist  normal complex spaces \ $\tilde{\tilde{M}}$ \ and  \ $\tilde{\tilde{S}}$, holomorphic maps  \ 
  $ g : \tilde{\tilde{M}}\rightarrow  \tilde{\tilde{S}}, \quad \tilde{\tilde{q}} :  \tilde{\tilde{S}} \to S, \quad \tilde{\tilde{\tau}} :  \tilde{\tilde{M}} \to M$ \ with the following properties :
  \begin{enumerate}[i)]
  \item The holomorphic map \ $g$ \ is geometrically f-flat with irreducible generic fibers\footnote{In the sense that there exists an open dense set where each fiber is irreducible (and reduced).}.
  \item The map \ $\tilde{\tilde{\tau}}$ \ is a proper modification.
    \item The map \ $\tilde{\tilde{q}}$ \ is proper, generically finite and surjective.
  \item We have \ $f\circ \tilde{\tilde{\tau}} = \tilde{\tilde{q}}\circ g$.
  \end{enumerate}
   \end{thm}
  
  \parag{Proof} Consider first a modification \ $\tau : \tilde{S} \to S$ \ with \ $\tilde{S}$ \ normal, such that the strict transform \ $\tilde{f} : \tilde{M} \to \tilde{S} $ \ is geometrically f-flat. Then let \ $\nu : \tilde{\tilde{M}} \to \tilde{M}$ \ be the normalisation, and apply the corollary \ref{utile} and the theorem \ref{Stein f-flat} to obtain the following commutative diagram
  $$ \xymatrix{ \tilde{\tilde{M}} \ar[d]^g \ar[dr]^{\tilde{\tilde{f}}} \ar[r]^{\nu} & \tilde{M} \ar[d]^{\tilde{f}} \ar[r]^{\tilde{\tau}} & M \ar[d]^f \\
  \tilde{\tilde{S}} \ar[r]^{\tilde{q}} & \tilde{S} \ar[r]^{\tau} & S} $$
  where \ $g$ \ is geometrically f-flat  with irreducible generic fiber, $\tilde{q}$ \ is proper finite and surjective, $\nu$ \ is the normalization map of \ $\tilde{M}$ \ and \ $\tilde{\tau}$ \ and \ $\tau$ \ are proper modifications. To conclude define \ $\tilde{\tilde{\tau}} : = \tilde{\tau}\circ \nu$ \ and \ $\tilde{\tilde{q}} : = \tau\circ \tilde{q}$. $\hfill \blacksquare$\\

\section{Appendix.}

The aim of this appendix is to show how the generalization by D. Mathieu [M.00]  of Khulmann's semi-proper direct image theorem to the case where the target space is a Banach analytic set has the following consequence :

\begin{thm}\label{semi-proper d.i. f-cycles}
Let \ $S$ \ and \ $M$ \ be reduced complex spaces, and denote by  \\ 
$\varphi : S \to \mathcal{C}_n^f(M)$ \ the classifying map of a f-analytic family of \ $n-$cycles in \ $M$. Assume that \ $\varphi$ \ is semi-proper. Then \ $\varphi(S)$ \ has a natural structure of weakly normal complex space such that the tautological family of \ $n-$cycles  parametrized by \ $\varphi(S)$ \  is an f-analytic family of cycles in \ $M$.
\end{thm}

Of course this result has an easy corollary which is the following property of  "universal reparametrization " (see the universal reparametrization theorem of [M.00]).

\begin{cor}\label{reparametrization}
In the situation of the previous theorem, consider an f-analytic family \ $(Y_t)_{t \in T}$ \ of \ $n-$cycles in \ $M$ \ parametrized by a weakly normal complex space \ $T$ \  such that for any \ $t \in T$ \ there exists an \ $s \in S$ \ such that \ $Y_t = X_s$. Then there exists an unique holomorphic map 
$$ \gamma : T \to \varphi(S) $$
such that we have \ $Y_t = X_{\gamma(t)}$ \ for all \ $t \in T$.
\end{cor}

\parag{Proof of the theorem} We know that  \ $\varphi(S)$ \ is closed and locally compact. Let \ $s_0 \in S$ \ and choose  relatively compact open set \ $W' \subset\subset W \subset\subset M$ \ meeting all irreducible components of \ $C_0 : = \varphi(s_0)$. Cover \ $\bar W$ \ by a finite set of adapted scales \ $E_i : = (U_i, B_i,j_i), i \in I$, and let \ $k_i : = \deg_{E_i}(C_0)$. Define the open neighbourhood 
$$ \mathcal{V} : = \Omega(W') \bigcap \big(\cap_{i \in I} \Omega_{k_i}(E_i)\big)$$
and consider the maps
$$\varphi' : \varphi^{-1}(\mathcal{V}) \to \mathcal{V}, \quad  \alpha : \mathcal{V} \to \prod_{i \in I} \ H(\bar U_i, \Sym^{k_i}(B_i)), \quad \psi : = \alpha\circ\varphi' .$$
Then \ $\varphi'$ \  is semi-proper as we know that  \ $\varphi(S)$ \ is  locally compact. It is now easy to check that  \ $\varphi(S')$ \  is homeomorphic to \ $\psi(S')$ \ by construction, using the characterization of compact sets in \ $\mathcal{C}_n^f(M)$. Now the semi-proper direct image of D. Mathieu [M.00]  gives a natural weakly normal complex structure on \ $\psi(S')$. $\hfill \blacksquare$\\

\parag{Remark} Using the proof above it is easy to show that, if we have an f-analytic family \ $(X_s)_{s \in S}$ \ with classifying map \ $\varphi : S \to \mathcal{C}_n^f(M)$ \ and an holomorphic map \ $g : S \to T$ \ between reduced complex spaces, such that \ $\varphi\times g$ \ is semi-proper, then there exists a natural structure of weakly normal complex space on \ $(\varphi\times g)(S)$ \ such the projection on \ $T$ \ is holomorphic and the projection on \ $\mathcal{C}_n^f(M)$ \ is the classifying map of an f-analytic family of cycles in \ $M$.\\

\parag{References}
\begin{enumerate}

 \item{[B.75]} Barlet, D. \textit{Espace analytique r\'eduit $\cdots$} S\'eminaire F.Norguet 1974-1975, L.N. 482 (1975), p.1-158.
 
 \item{[B.78]} Barlet, D. \textit{ Majoration du volume $\cdots$ } S\'eminaire Lelong-Skoda 1978-1979, L.N. 822 (1980), p.1-17.
 
 \item{[B.08]} Barlet, D. \textit{Reparam\'etrisation universelle de familles f-analytiques de cycles et f-aplatissement g\'eom\'etrique} Comment. Math. Helv. 83 (2008),\\ p. 869-888.
  
 \item{[BOOK]} Barlet, D. and  Magnusson, J. \textit{Book to appear}.
 
 \item{[B.Mg.10]} Barlet, D. and  Magnusson, J. article  in preparation.
 
 \item{[Bi.64]} Bishop, E \textit{Conditions for the analyticity of certain sets} Mich. Math. J. 11 (1964), p. 289-304.
 
  \item{[C.60]} Cartan, H. \textit{Quotients of complex analytic spaces} International Colloquium on Function Theory, Tata Institute (1960), p.1-15.
  
  \item{[D.90]} Dethloff, G. \textit{A new proof of a theorem of Grauert and Remmert by\ $ L_2$ \ methods}, Math. Ann. 286 (1990) p.129-142.
 
 \item{[G.83]} Grauert, H. \textit{Set theoretic complex equivalence relations} Math. Ann. 265 (1983), p.137-148.
 
 \item{[G.86]} Grauert, H. \textit{On meromorphic equivalence relations} Proc. Conf. Complex Analysis, Notre-Dame (1984) Aspects Math. E9 (1986), p.115-147.
 
 \item{[G.R.58]} Grauert, H. and Remmert, R. \textit{Komplexe R{\"a}ume}, Math. Ann. 136 (1958) p. 245-318.
 
 \item{[K.64]} Kuhlmann, N. \textit{{\"U}ber holomorphe Abbildungen komplexer R{\"a}ume} Archiv der Math. 15 (1964), p.81-90.
 
 \item{[K.66]} Kuhlmann, N. \textit{Bemerkungen {\"u}ber  holomorphe Abbildungen komplexer R{\"a}ume} Wiss. Abh. Arbeitsgemeinschaft Nordrhein-Westfalen 33, Festschr. Ged{\"a}achtnisfeier K. Weierstrass (1966), p.495-522.
 
 \item{[M.00]} Mathieu, D. \textit{Universal reparametrization of a family of cycles : a new approach to meromorphic equivalence relations}, Ann. Inst. Fourier (Grenoble) t. 50, fasc.4 (2000) p.1155-1189.
  
 \item{[R.57]} Remmert, R. \textit{Holomorphe und meromorphe Abbildungen komplexer R{\"a}ume}, Math. Ann. 133 (1957), p.328-370.
 
 \item{[S.93]} Siebert, B. \textit{Fiber cycles of holomorphic maps I. Local flattening} Math. Ann. 296 (1993), p.269-283. 
 
  \item{[Si.94]} Siebert, B.  \textit{Fiber cycle space and canonical flattening II} Math. Ann. 300 (1994), p. 243-271.

\end{enumerate}

\end{document}